\documentclass[12pt]{amsart}
\usepackage{amsfonts}
\usepackage{amssymb}
\usepackage{graphicx}
\usepackage[usenames]{color}
\usepackage{soul}

\def\D{{\mathbb D}}  
\def\C{{\mathbb C}}  
\def\R{{\mathbb R}}

\def\({\left(}       \def\){\right)}

\newtheorem{prop}{\sc Proposition}

\newtheorem{thm}{\sc Theorem}
\newtheorem{cor}{\sc Corollary}
\newtheorem{ex}{\sc Example}

\newtheorem{other}{\sc Theorem}              
\newenvironment{pf}{\noindent{\textit{Proof. }}}{$\Box$ }

\begin{document}
\title[Always-convex harmonic shears]
{Always-convex harmonic shears}

\author[ ]{Rodrigo Hern\'andez}
\address{Facultad de Ingenier\'{\i}a y Ciencias, Universidad Adolfo Ib\'a\~nez. Av. Padre Hurtado 750, Vi\~na del Mar, Chile.} \email{rodrigo.hernandez@uai.cl}
\author[ ]{Mar\'{\i}a J. Mart\'{\i}n}
\address{Departamento de An\'alisis Matemático \& IMAULL, Universidad de La Laguna. Av.  Astrofísico Francisco Sánchez, s/n. Facultad de Matemáticas. 38200, La Laguna, Tenerife, Spain.} \email{maria.martin@ull.es}
\author[ ]{Fernando~P\'erez-Gonz\'alez}
\address{Departamento de An\'alisis Matemático \& IMAULL,  Universidad de La Laguna. Av.  Astrofísico Francisco Sánchez, s/n. Facultad de Matemáticas. 38200, La Laguna, Tenerife, Spain.} \email{fpergon@ull.edu.es}
\author[ ]{Magdalena  Wo{\l}oszkiewicz-Cyll}
\address{Instytut Informatyki i Matematyki, Uniwersytet Marii Curie-Sk{\l}odowskiej. Pl. M. Curie-Sk{\l}odowskiej 1,
20-031 Lublin, Poland.} \email{magdalena.woloszkiewicz-cyll@mail.umcs.pl}

\date{\today}

\thanks{R. Hern\'andez and M. J. Mart\'in were partially supported by Spanish  MICIU research projectPID2024-160185NB-I00 MICIU research project PID2024-160185NB-I00.; M. Wo{\l}oszkiewicz-Cyll was supported by a grant of National Science Center, Poland (NCN), 2022/06/X/ST1/01611.}

\begin{abstract} We determine completely the analytic functions $\varphi$
in the unit disk $\mathbb D$ such that for all (normalized)
orientation-preserv\-ing harmonic mappings $f=h+\overline g$ produced by the shear construction with
$h+ g=\varphi$,  the condition that each $f$ maps $\mathbb D$ onto a
convex domain holds. As a consequence, we obtain the following more general result: for a given complex number $\eta$, with $|\eta|=1$, we characterize those holomorphic mappings $\varphi$ in $\mathbb D$ such that every harmonic function $f=h+\overline g$ as above with $h-\eta g=\varphi$ maps $\mathbb D$ onto a convex domain. 
The resulting functions are mappings onto a half-plane and mappings onto a strip, and the shear direction, determined by the parameter $\eta$ above, is parallel to the linear boundaries of the half-planes and strips.

\end{abstract}

\maketitle
\date{\today}
\section*{Introduction}

A domain $\Omega$ in the complex plane is said to be \emph{convex in the $t$ direction}, $t\in[0, \pi)$, if the intersection of $\Omega$ with any line parallel to the line through $0$ and $e^{it}$ is either an interval or empty. A convex domain is a domain convex in every direction. The univalent complex-valued function $f$ in the unit disk $\D$ is \emph{convex in the $t$ direction} (\emph{resp.} \emph{convex}) if $f(\D)$ is convex in the $t$ direction (\emph{resp.} convex). Every domain convex in some direction is a close-to-convex domain (see, for instance, \cite{Pom-CtC} and the references therein).
\par\smallskip
Given any holomorphic function $f$ in the unit disk with $f(0)=1-f'(0)=0$ and such that $f(\D)$ is convex in the direction $t$, the rotation $f_\gamma(z)=e^{-\gamma i} f(e^{\gamma i}z)$, with $\gamma=t-\pi/2$, maps $\D$ onto a domain convex in the vertical direction (that is, in the $\pi/2$ direction)  and satisfies the same normalizations. Hence, we can focus on this latter family of analytic mappings, denoted by $CV(i)$.

\par\smallskip
Fejer \cite[p.~61]{F} was the first to call attention to this family, though Robertson \cite{R1, R2} was the first to make an intensive study of $CV(i)$. According to \cite[p.~196]{G}, it was a great surprise to the mathematical community to learn that $f\in CV(i)$ does not imply that the dilation $f_r$, $0<r<1$, defined by $f_r(z)=f(rz)/r$, is also a member of this class \cite{HS}. An explicit example of a function that is convex in the vertical direction but whose restriction to the disk $|z|<r$ does not have this property for any radius $\sqrt 2-1<r<1$ was constructed by Goodman and Saff \cite{GS}. They conjectured that the radius $\sqrt 2-1$ is best possible, which was later proved by Ruscheweyh and Salinas \cite{RS}.
\par\smallskip
A complex-valued function $f$ is harmonic if both the real and imaginary parts of $f$ are harmonic. Every harmonic mapping $f$ in the unit disk has a \emph{canonical decomposition} $f=h+\overline g$, where $h$ and $g$ are analytic in $\mathbb{D}$ with $g(0)=0$. A locally univalent function $f$  is called  \emph{orientation-preserving} if its Jacobian $J_f=|h'|^2-|g'|^2$ is positive in $\mathbb{D}$  \cite{Lewy}, so that    $\omega=g'/h'$  is analytic in $\mathbb{D}$ and  $\omega(\mathbb{D})\subset\mathbb{D}$. If $J_f<0$ in $\mathbb{D}$, $f$ is   \emph{orientation-reversing}. It is obvious that $f$ is orientation-reversing if and only if $\overline f$ is orientation-preserving.

\par\smallskip
 Let $S_H$ be the class which consists of those orientation-preserving univalent harmonic mappings $f=h+\overline g$ in the unit disk normalized by the conditions $h(0)=g(0)=0$ and $h'(0)=1$. The class $S_H$ is known to be normal whereas the class $S_H^0=\{f=h+\overline g\in S_H\colon g'(0)=0\}$ is compact (see \cite[Ch.~5]{D-Harm}). The class $\mathcal S$ of analytic and univalent functions $\varphi$ in the unit disk normalized by the conditions $\varphi(0)=1-\varphi'(0)=0$ is a subfamily of $S_H^0$.  Notice that the requirement on the normalizations $h(0)=g(0)=1-h'(0)=g'(0)=0$ can be always obtained by applying appropriate (invertible) affine harmonic transformations as described in \cite[Ch.~5]{D-Harm}.

\par\smallskip
The dilatation $\omega=g'/h'$ of any $f\in S_H^0$ is a holomorphic function in the unit disk that fixes the origin and satisfies $\omega(\mathbb{D})\subset \mathbb{D}$. We will say that any holomorphic function $\omega$ in $\mathbb{D}$ satisfying these two properties is a \emph{Schwarz function}.

\par\smallskip
The shear construction, introduced by Clunie and Sheil-Small \cite{CS-S}, is an effective method for producing univalent harmonic mappings in $\D$ with a specified dilatation onto domains convex in one direction by ``shearing'' a univalent analytic mapping along parallel lines. It can be easily verified that the version presented here is equivalent to the one in \cite[Thm.~5.3]{CS-S}.
\begin{other}\label{thm-shear}
Let $f=h+\overline g$ be an orientation-preserving harmonic mapping in the unit disk. Then $f$ is univalent and maps the unit disk onto a domain convex in the vertical direction (that is, in the direction $t=\pi/2$) if and only if the analytic function $h+g$ maps the unit disk onto a domain convex in the vertical direction.
\end{other}
\par
A direct application of Theorem~\ref{thm-shear} shows that given any $\varphi \in CV(i)$ and any Schwarz function $\omega$, the solution $(h,g),$ determined by the linear system of equations
\begin{equation}\label{eq-system}
\left\{\begin{array}{cc}
h+g=\varphi\\
g^\prime/h^\prime= \omega
\end{array}
\right.\,,\quad h(0)=g(0)=0\,,
\end{equation}
gives rise to the harmonic mapping $f=h+\overline g \in S_H^0$ which maps $\D$ onto a domain convex in the vertical direction and has dilatation $\omega$.

\par\smallskip
In other words, it could be said that the shears of a function $\varphi\in CV(i)$ are \emph{always convex in the vertical direction} (and univalent) since every solution $f=h+\overline g$ which satisfies \eqref{eq-system} is convex in the vertical direction \emph{for every} Schwarz function $\omega$.

\par\smallskip
It is easy to check that, in terms of $\varphi$ and $\omega$, the function $f=h+\overline g$ obtained from \eqref{eq-system} has the form
\[
f(z)={\rm Re\,} \varphi(z)+i\, {\rm Im\,} \int_0^z \varphi'(\zeta) \frac{1-\omega(\zeta)}{1+\omega(\zeta)}\, d\zeta\,.
\]

\par\medskip
It is known that every orientation-preserving harmonic mapping $f$ which maps the unit disk onto a convex domain is univalent \cite[Thm.~5.7]{CS-S}. The class $C_H^0$ consists of harmonic mappings $f=h+\overline g \in S_H^0$ which map $\D$ onto a convex domain.

\par\smallskip
The following theorem was proved in  \cite{CS-S}. We state the result for future reference.

\begin{other}\label{thm-convexcharact}
Let $f=h+\overline g\in S_H^0$. Then $f$  maps the unit disk onto a convex domain if and only if for each $t\in[0,\pi)$ the analytic function $h-e^{2it}g$  is convex in the direction $t$.  
\end{other}   

\par\smallskip
It is worth noting that the condition characterizing convexness, as stated in the theorem, is not easy to verify in practice. Nevertheless, Hengartner and Schober in \cite{HS2} managed to show that it is satisfied by any (normalized) orientation-preserving harmonic mapping $f=h+\overline g$ obtained from the solution to \eqref{eq-system} with
\[
\varphi(z)=\frac{-i}{2}\log\frac{1+iz}{1-iz}
\,,\quad z\in\D\,,
\] 
and $f$ therefore maps $\D$ onto a convex domain \emph{for every} Schwarz function $\omega$.
\par\smallskip
Boyd \emph{et al.} \cite{B}, adapted the arguments in \cite{HS2} to prove that if $f=h+\overline g$ is obtained from \eqref{eq-system} with $\varphi(z)=z/(1-z)$ for $z\in\D$,
then, again, $f$ maps the unit disk onto a convex domain \emph{for every} Schwarz function $\omega$.

\par\smallskip
The effectiveness of the shear construction in producing $C_H^0$ functions is, however, limited. A simple example that justifies this assertion is provided by the choices $\varphi(z)=z$ and $\omega(z)=-iz$ in \eqref{eq-system}. The resulting function $f$ satisfies that $f(z)=-iF(iz)$, $z\in\D$, where $F$ is the function in the example in \cite[p. 38]{D-Harm}. The function $F$ has the range shown in \cite[Figure 3.1, p. 39]{D-Harm}, which is not convex. Therefore, $f(\D)$ is not a convex domain.

\par\smallskip
More generally, we can note that a given function $f=h+\overline g\in S_H^0$ is completely determined by the system of equations
\begin{equation}\label{eq-system0}
\left\{\begin{array}{cc}
h-\eta g=\varphi\\
g^\prime/h^\prime= \omega
\end{array}
\right.\,,\quad h(0)=g(0)=0\,,
\end{equation}
where $\eta=e^{2i\theta}\in \partial \D$ and $0<2\theta\leq 2\pi$.  We call such function $f$ \emph{the shear of $\varphi$ in the $\theta$ direction with dilatation $\omega$}. It is evident that different dilatations $\omega$ in \eqref{eq-system0} give rise to different shears of $\varphi$ in the $\theta$ direction. We say that the shears of $\varphi$ in the $\theta$ direction are \emph{always-convex} if the harmonic mappings $f=h+\overline g$ obtained from \eqref{eq-system0} with dilatation $\omega$ are convex \emph{for all} Schwarz functions $\omega$.  Within these terms, the results in \cite{B} mentioned above can be stated as follows: the shears of 
\begin{equation}\label{eq-H}
\mathcal H(z)=z/(1-z)\,,\quad  z\in\D\,, 
\end{equation}
in the $\pi/2$ direction are always-convex. 
\par\smallskip

By applying the same arguments as in \cite{B} it was proved in \cite{FHM} that vertical shears of the functions
\begin{equation}\label{eq-llambda}
\mathcal L_\lambda(z)=\frac{1}{2i{\rm Im} \lambda} \log\frac{1-\overline\lambda z}{1-\lambda z}\,,\quad z\in\D\,,
\end{equation}
where $\lambda\in\partial\D\setminus\{-1,1\}$ are always-convex as well.
\par\smallskip
It should  be noted that there exists a connection between the parameter $\theta$ and the function $\varphi$ in \eqref{eq-system0} in order to obtain always-convex harmonic shears of the function $\varphi$ in the $\theta$ direction: even though the shears of $\mathcal H$ in the $\pi/2$ direction are always-convex, the shears of $\mathcal H$ in the $0$ direction are not,  as the following fundamental example shows. We refer the reader to \cite{4} for the details on the relevance on the function we now present. It is worth mentioning also the nice articles \cite{2} and \cite{3}, where the authors obtain different important consequences by using the strengthness of Theorem~\ref{thm-shear} (among other results from the literature combined in an appropriate way).

\par\smallskip

\begin{ex}\label{ex-f_0}
 Let $f_0=h_0+\overline{g_0}\in S_H^0$ be defined, for $z$ in $\D$, by the solution to 
\[
\left\{\begin{array}{cc}
h_0-g_0=\mathcal H\\
g_0^\prime/h_0^\prime= I
\end{array}
\right.\,,\quad h(0)=g(0)=0\,,
\]
where $I(z)=z$ for $z$ in the unit disk and $\mathcal{H}$ is as in \eqref{eq-H}.

In \cite{4}, it is proved that $f_0$ is not starlike and hence, it is not convex (see \cite[Fig. 7]{4} for a representation of the image of  $f_0$, which is denoted by $f_3$ in the article mentioned).

\par\smallskip

\end{ex}

This last example shows that there is a Schwarz function ($\omega=I$) for which the horizontal shear of $\mathcal H$ with dilatation $\omega$ does not map the unit disk onto a convex domain. Hence, as announced before, the shears of $\mathcal H$ in the $0$ direction are not always-convex.
   
\par\medskip

The main goal in this paper is to characterize \emph{all} pairs $(\eta, \varphi)$, where $\eta=e^{2i\theta} $,  such that the shears of $\varphi$ in the $\theta$ direction are always-convex. 

\par\smallskip

It should be noted that, since the function $\omega\equiv 0$ is a Schwarz function (\emph{i.e.}, $\omega$ is analytic in $\D$, $\omega(0)=0$, and $\omega(\D)\subset \D$), if we assume that the shears of $\varphi$ in the $\theta$ direction are always-convex, then, in particular, the datum $\varphi$ in \eqref{eq-system0} must be convex. This property will provide fundamental tools to prove our main result, Theorem~\ref{thm-main} below. More precisely, different theorems from other authors will be required in the proof of Theorem~\ref{thm-main}. We mention some of the results due to Brannan \cite{Brannan}, Noonan \cite{Noonan}, and Pinchuk \cite{Pinchuk}, from the theory of analytic functions of bounded boundary rotation; some theorems due to Dorff \cite{Dorff} and Grygorian and Nowak \cite{GN}, on the structure of some particular convex harmonic mappings; and different theorems by Abu-Muhanna and Schober \cite{AM-S}, Bshouty, Lyzzaik, and Weitsman \cite{BLW}, Clunie and Sheil-Small \cite{CS-S}, Hengartner and Schober \cite{HS3},  Laugesen \cite{Laugesen}, and  Sheil-Small \cite{S-S}, on the boundary behaviour of harmonic mappings in the plane. The fact that $\varphi$ is assumed to map the unit disk conformally onto a convex domain, jointly with some other auxiliary lemmas and theorems that will be proved along the paper, will allow us to combine the results mentioned and characterize the shears of $\varphi$ in the $\theta$ direction that are always-convex.  

\section{The use of rotations}

As mentioned above, our main goal is to characterize the pairs $(\eta,\varphi)$, where $\eta=e^{2i\theta}$, $0< 2\theta\leq 2\pi$, has modulus one and $\varphi\in\mathcal S$ such that any harmonic mapping $f=h+\overline g$ obtained from the solution to the system 
\begin{equation}\label{eq-system1}
\left\{\begin{array}{cc}
h-\eta g=\varphi\\
g^\prime/h^\prime= \omega
\end{array}
\right.\,,\quad h(0)=g(0)=0\,,
\end{equation}
maps the unit disk onto a convex domain for every Schwarz function $\omega$ in \eqref{eq-system1}. In the cases when this property is satisfied, we say that the shears of $\varphi$ in the direction $\theta$ are always-convex and we call the pair $(\eta, \varphi)$ \emph{admissible}. As mentioned above, the function $\varphi$ in any admissible pair $(\eta, \varphi)$ must be convex analytic since $\omega\equiv 0$ is a Schwarz function.

\par\smallskip

Several examples of admissible pairs of the form $(-1, \varphi)$ have been given above. Other examples of admissible pairs $(\eta, \varphi)$ for $\eta$ not necessarily equal to $-1$ are provided by the next proposition. Here, for $\xi\in\partial \D$ and $f \in S_H^0$ we define the \emph{rotation}
\begin{equation}\label{eq-rotation}
f_\xi(z)=\overline \xi f(\xi z)\,,\quad z\in\D\,.
\end{equation}
It is obvious that $(f_\xi)_{\overline \xi}=(f_{\overline{\xi}})_{\xi} =f$.
\begin{prop}\label{prop-rotations}
Let $\eta$ be a complex number of modulus one and let $\varphi$ be a convex analytic function (and hence necessarily univalent) in $\D$ with $\varphi(0)=1-\varphi'(0)=0$. Then, the pair $(\eta, \varphi)$ is admissible if and only if the pair $(\eta\xi^2, \varphi_{\overline\xi})$ is admissible for each $|\xi|=1$. 
\end{prop}
\par\smallskip
\begin{pf}
Suppose that the pair $(\eta, \varphi)$ is admissible, so that the solution $f=h+\overline g$ to the system \eqref{eq-system1} turns out to be convex for every Schwarz function $\omega$. In order to get a contradiction, assume that there exist $\xi\in\partial\D$ and a Schwarz function $\omega_0$ so that the solution $F=H+\overline G$ to
\[
\left\{\begin{array}{cc}
H-\eta\xi^2 G=\varphi_{\overline\xi}\\
G^\prime/H^\prime= \omega_0
\end{array}
\right.\,,\quad H(0)=G(0)=0\,,
\]
is  not convex.
\par\smallskip
Consider the rotation $F_\xi$ defined by \eqref{eq-rotation}. The canonical decomposition of $F_\xi=H^\xi+\overline {G^\xi}$ is given by
\begin{equation*}\label{eq-hyglambda}
H^\xi(z)=H_\xi(z)=\overline\xi H(\xi z)\quad\text{and}\quad G^\xi(z)=\xi G(\xi z)\,,\quad z\in\D\,.
\end{equation*}
Since $F$ is not convex, the function $F_\xi$ is not convex either. However, a straightforward calculation shows that ($H^\xi, G^\xi$) solve \eqref{eq-system0} with datum $\omega$ defined in $\D$ by $\omega(z)=\xi^2\omega_0(\xi z)$. This gives the desired contradiction and shows that if $(\eta, \varphi)$ is admissible then $(\eta\xi^2, \varphi_{\overline\xi})$ is also admissible. It is obvious that the reverse implication also holds since the transformation $(\eta, \varphi) \mapsto (\eta\xi^2, \varphi_{\overline\xi})$ can be inverted.
\end{pf}
\par\smallskip

The relevance of Proposition~\ref{prop-rotations} for our purposes is that we can determine the admissible pairs $(\eta, \varphi)$ as rotations of the admissible pair $(-1,\varphi)$. In other words, we have that $(-1,\varphi)$ is admissible if and only if $(\eta, \varphi_{\overline\xi})$, where $\xi^2=-\eta$, is. Hence, the problem to be solved is the determination (up to rotations) of all those locally univalent analytic functions $\varphi$ in $\D$ such that the harmonic mapping $f=h+\overline g$, where $h$ and $g$ solve the system
\begin{equation}\label{eq-systemimportant}
\left\{\begin{array}{cc}
h+ g=\varphi\\
g^\prime/h^\prime= \omega
\end{array}
\right.\,,\quad h(0)=g(0)=0\,,
\end{equation}
is convex independently of the Schwarz function $\omega$.
\par\smallskip
We again point out that if $(-1, \varphi)$  is admissible, then $\varphi$ must map the unit disk onto a convex domain, as explained above. 

\par\smallskip
Consider the functions $\mathcal H$ and $\mathcal L_\lambda$ defined by \eqref{eq-H} and \eqref{eq-llambda}, respectively, together with the rotation $\mathcal H_{-1}$ as in \eqref{eq-rotation} of the function $\mathcal  H$, in order to define the set 
\begin{equation}\label{eq-classA}
\mathcal A=\left\{\mathcal H, \mathcal H_{-1}, \mathcal L_\lambda, \text{ where }\lambda\in\partial\D\setminus\{-1,1\}  \right\}\,.
\end{equation}
Notice that $(\mathcal L_\lambda)_{-1}=\mathcal L_{-\lambda}$  and the class $\mathcal A$ is therefore closed under rotations of the form \eqref{eq-rotation} for $\xi=-1$.
\par\smallskip

A direct use of Proposition~\ref{prop-rotations} also shows that the pair $(-1, \varphi)$ is admissible if and only if the pair $(-1, \varphi_{-1})$ is admissible, where $\varphi_{-1}$ is the rotation defined as in \eqref{eq-rotation}. This remark, along with the results in \cite{B} and \cite{FHM}, provides a direct proof of the following result. We omit the details.

\begin{prop}\label{prop-model}
    Let $\varphi\in\mathcal A$. Then $(-1, \varphi)$ is admissible. That is, the vertical shears of $\varphi$  are always-convex. 
\end{prop}

\par\smallskip

A natural question that arises at this point is whether other rotations $\varphi_\xi$, where $\varphi\in\mathcal A$ and $\xi\in\partial \D$, $\xi\neq-1$, give rise to admissible pairs $(-1, \varphi_\xi)$. 

\par\smallskip
In order to answer this question, it is worth noting that if a given harmonic mapping $f=h+\overline g \in S_H^0$ has the property that $h-g$ is not convex in the vertical direction, then the harmonic function $h-\overline g$ cannot be convex. An easy proof of this assertion follows from the fact that if $h-\overline g$ is convex, then $h+e^{2t i} g$ must be convex in the $t$ direction for each $t\in[0,\pi)$, according to Theorem~\ref{thm-convexcharact}. This property does not hold for $t=\pi/2$ unless $h-g \in CV(i)$. This simple observation is the key tool to prove our next result.

\begin{thm}\label{thm-chd}
Let $\varphi\in\mathcal A$ and let $\xi$ be a complex number of modulus $1$. The vertical shears $f=h+\overline g$ of the rotation $\varphi_\xi$ of $\varphi$, as defined in \eqref{eq-rotation}, are always-convex if and only if $\xi\in\{-1,1\}$.
\end{thm}
\begin{pf}
The sufficiency of the condition in the statement follows from Proposition~\ref{prop-model}.

\par\smallskip
To prove the necessity we first notice that since vertical shears of $\varphi$ are always-convex if and only if vertical shears of $\varphi_{-1}$ are, we can assume without loss of generality that either $\varphi=\mathcal H$ or $\varphi=\mathcal L_\lambda$ for some $\lambda\in\partial{\D}$ with a positive imaginary part. In fact, using again that by Proposition~\ref{prop-rotations}, we have that $(-1, \varphi_{-\xi})$ is admissible if and only if $(-1, \varphi_\xi)$ is admissible, we can assume that $\xi=e^{i\theta}$ for some $\theta\in [0, \pi)$. Thus, we are to prove that $\xi=1$.
\par\smallskip
Let us first consider the case $\varphi=\mathcal H$. In order to get the result for vertical shears of rotations of the function $\mathcal H$, let us assume that whenever $\omega$ is a Schwarz function and $(h, g)$ solves
\[
\left\{\begin{array}{cc}
h+g=\mathcal H_{\xi}\\
g^\prime/h^\prime= \omega
\end{array}
\right.\,,\quad h(0)=g(0)=0\,,
\]
the harmonic mapping $f=h+\overline g$ maps the unit disk onto a convex domain. In particular, we have that the function $F=H+\overline G$ that is produced by the solution ($H$, $G$) to the system 
\[
\left\{\begin{array}{cc}
H+G=\mathcal H_{\xi}\\
G^\prime/H^\prime= -\xi I
\end{array}
\right.\,,\quad H(0)=G(0)=0\,,
\]
must be convex.

\par\smallskip
Now, recall that the harmonic mapping $f_0=h_0+\overline{g_0}$ defined in Example~\ref{ex-f_0} is not convex. In addition, the functions $h_0$ and $g_0$ in the canonical decomposition of $f_0$ satisfy
\[
\left\{\begin{array}{cc}
h_0-g_0=\mathcal H\\
g_0^\prime/h_0^\prime= I
\end{array}
\right.\,,\quad h_0(0)=g_0(0)=0\,.
\]

\par\smallskip
A straightforward computation shows
\[
h_0'(\xi z)(1-\xi z)=\frac{1}{(1-\xi z)^2}=\mathcal H_{\xi}'(z)= H'(z)\,(1-\xi z)\,, \quad z\in\D\,.
\]
Thus, for all $z$ in the unit disk, $H'(z)=h_0'(\xi z)$ and $H(z)=\overline \xi h_0(\xi z)$. Also,
\[
-\xi z= \frac{G'(z)}{H'(z)}=\frac{G'(z)}{h_0'(\xi z)}= \frac{g_0'(\xi z)}{h_0'(\xi z)}\, \frac{G'(z)}{g_0'(\xi z)}=\xi z\, \frac{G'(z)}{g_0'(\xi z)}\,,\quad z\neq 0\,.
\]
This means that for all such $z$ (hence, for all $z\in\D$ by the identity principle for analytic functions), $G'(z)=-g_0'(\xi z)$, which gives $G(z)=-\overline\xi g_0(\xi z)$.
\par\smallskip
In other words, the convex function $F$ is given by
\[
F(z)=H(z)+\overline{G(z)}= \overline \xi h_0(\xi z) -\xi \overline{ g_0(\xi z)}\,, \quad z\in\D\,.
\]
Since $F$ is a convex function, its rotation
\[
F_{\overline \xi}(z)= h_0(z) - \xi^2\overline{g_0(z)}\,, \quad z\in\D\,,
\]
must be convex as well. Recall that $\xi=e^{i\theta}$ for some $\theta\in [0, \pi)$ and, by \cite[Thm. 5.7]{CS-S}, the harmonic mapping $h_0-\xi^2\overline g_0$ is convex if and only if for all $t\in[0, \pi]$ the analytic functions $h_0+e^{2it}\overline\xi^2\, g_0$ are convex in the $t$ direction. This means that, in particular, for $t=\theta$, the function $h_0+ g_0$ must be convex in the direction $\theta$. But $h_0+g_0=k$, where $k$ is the \emph{Koebe function} defined in the unit disk by $k(z)=z/(1-z)^2$. The domain $k(\D)$ is  convex only when $\theta=0$ and hence $\xi=1$.
\par\smallskip
In cases where $\varphi=\mathcal L_\lambda$ for some $|\lambda|=1$ with a positive imaginary part, we can repeat the same argument used for the case $\varphi=\mathcal H$ above. The statement $\xi=1$ will be proved once we prove that for all such $\lambda$, the vertical shear $f=h+\overline g$ of $\mathcal L_\lambda$ with dilatation $-I$ satisfies that $h-g$ is convex only in the horizontal direction.
\par
The functions $h$ and $g$ of these vertical shears satisfy
\[
\left\{\begin{array}{cc}
h+g=\mathcal L_{\lambda}\\
g^\prime/h^\prime= -I
\end{array}
\right.\,,\quad h(0)=g(0)=0\,.
\]
Hence,
\[
h'(z)(1-z)=h'(z)+g'(z)=\frac{1}{(1-\lambda z)(1-\overline\lambda z)}\,.
\]
This gives
\[
h'(z)-g'(z)=h'(z) (1+z)=\frac{1+z}{1-z} \frac{1}{(1-\lambda z)(1-\overline \lambda z)}\,.
\]
The result is now evident from the Schwarz-Christoffel formula and the fact that for all $z\in\D$, $(h-g)'(z)=\overline{(h-g)'(\overline z)}$. However, we have decided to include for the sake of completeness an alternative way of concluding that $h-g$ is convex only in the horizontal direction.
\par
Integrating the previous equation and using the fact that $h(0)=g(0)=0$, we get
\[
(h-g)(z)=\frac{1}{2-2{\rm Re }\lambda} \log\frac{(1-\lambda z)(1-\overline\lambda z)}{(1-z)^2}\,,
\]
where the principal branch of the argument has been chosen in order to define the logarithm.
\par
Now, notice that
\[
\frac{\partial}{\partial z} \left(\frac{(1-\lambda z)(1-\overline\lambda z)}{(1-z)^2}\right)= (2-2{\rm Re }\lambda)k'(z)\,,
\]
where, once more, $k$ is the Koebe function defined above.
Hence,
\[
\frac{(1-\lambda z)(1-\overline\lambda z)}{(1-z)^2}=(2-2{\rm Re }\lambda)k(z)+1\,,
\]
which transforms the unit disk onto the entire complex plane minus the slit on the real axis that contains the interval $(-\infty,(1+{\rm Re}\lambda)/2]$. Note that $(1+{\rm Re}\lambda)/2$ is a positive real number and that the principal branch of the logarithm maps the whole complex plane minus the negative part of the real axis conformally onto the strip $\{z\in\C:\,|{\rm Im }z| < \pi\}$ and the positive part of the real axis to the real axis. Hence, the mapping $h-g$ transforms the unit disk onto a strip with a horizontal slit from $-\infty$ to $\log((1+{\rm Re}\lambda)/2)$ and is therefore convex only in the horizontal direction.
\end{pf}
\par\smallskip
An important consequence of the previous proposition is summarized in the following corollary.
\begin{cor}\label{cor-consequence}
Let $A$ and $B$ be two complex numbers of modulus one. Assume that the analytic function $\psi$ in the unit disk satisfies
\[
\psi'(z)=\frac{1}{(1-Az)(1-Bz)}\,,\quad z\in \D.
\]
Then, the vertical shears of $\psi$ are always-convex if and only if $A=\overline B$.
\end{cor}
\begin{pf}
A straightforward calculation shows that if $A=\overline B$, then $\psi$ belongs to the set $\mathcal A$ given in \eqref{eq-classA}. Therefore, by Proposition~\ref{prop-model}, the vertical shears of $\psi$ are always-convex. 

\par\smallskip
To prove the sufficiency of the condition, let us suppose that such vertical shears are always-convex. We can then find $\lambda$ and $\xi$ of modulus one which solve the system
\[
\left\{\begin{array}{cc}
A=\lambda \xi\\
B=\overline\lambda \xi
\end{array}\right.
\]
to have that the vertical shears of the function $\psi$ with
\[
\psi'(z)=\frac{1}{(1-\lambda \xi z)(1-\overline\lambda \xi z)}\,,\quad z\in \D\,,
\]
are always-convex. 
\par\medskip
The previous equation shows that the function $\psi$ is equal to the rotation $\varphi_\xi$ of a function $\varphi\in\mathcal A$, as defined in \eqref{eq-rotation}. By Theorem~\ref{thm-chd}, $\xi\in\{-1, 1\}$ and therefore, $A=\overline B$.
\end{pf}

\section{Some properties of certain families of  analytic functions}
Bearing in mind that the datum $\varphi\in\mathcal S$ in any admissible pair $(\eta, \varphi)$ must be convex, it is convenient to review some of the properties of convex analytic functions that will be used below. 

\subsection{Convex analytic functions onto half-planes or strips} \label{sec-convexanalytic}

\par\smallskip
Any  function $\varphi\in\mathcal S$ that maps the unit disk onto a half-plane must be a M\"obius transformation. The conditions $\varphi(0)=1-\varphi'(0)=0$, jointly with the fact that $\varphi$ is analytic in $\D$ and maps the unit circle onto a straight line give that, for $z\in\D$,
\[
\varphi(z)=\frac{z}{1-cz}\,,\quad  |c|=1\,.
\]

The derivative of such a function is equal to
\[
\varphi'(z)=\frac{1}{(1-cz)^2}\,,\quad  |c|=1\,.
\]
Therefore, by Corollary~\ref{cor-consequence}, if the vertical shears of $\varphi$ are always-convex and $\varphi$ maps the unit disk onto a half-plane, then $\varphi\in\mathcal A$. 

\par\smallskip
A similar result is obtained if we assume that the vertical shears of $\varphi\in\mathcal S$ are always-convex and $\varphi$ maps the unit disk onto a strip. Namely, let $\varphi\in\mathcal S$ map the unit disk onto a strip. Bearing in mind that $\varphi(0)=0$, we can apply a rotation to get that, for some $\xi$ of modulus one, the function $\varphi_\xi$ defined in $\D$ as $\varphi_\xi(z)=\overline\xi\varphi(\xi z)$ maps $\D$ onto a strip of the form
\[
\mathbb S=\{z\in\C\colon a<{\rm Re\,} z <b\}\,, 
\]
where $a$ and $b$ are real numbers with $a<0<b$ and $|a|<b$.

\par\smallskip
It is an easy exercise to show that 
\[
\varphi_\xi(z)=K\log{\frac{1-\overline\lambda z}{1-\lambda z}}\,,\quad z\in\D\,,
\]
for some $K\in\C$ and some $\lambda$ of modulus one and different from $1$ and $-1$. This readily implies that 
\[
\varphi(z)=\xi K\log{\frac{1-\overline\lambda\,\overline\xi z}{1-\lambda\overline \xi z}}\,,\quad z\in\D\,,
\]
where $K=1/(2i{\rm Im\,} \lambda)$, since $\varphi'(0)=1$. Hence,
\begin{equation}\label{eq-strip}
\varphi'(z)=\frac{1}{(1-\overline\lambda\,\overline\xi z)(1-\lambda\,\overline\xi z)}\,.
\end{equation}

\par\smallskip

Another application of Corollary~\ref{cor-consequence} is that if the vertical shears of $\varphi$ are always-convex and $\varphi$ maps the unit disk onto a strip, then, necessarily, the parameter $\xi$ in \eqref{eq-strip} is either $-1$ or $1$, which implies that $\varphi\in\mathcal A$.

\par\smallskip
The previous examples are just a sample of the following well-known result proved, for instance, in \cite[Thm. 2]{Pom-CtC}: if an analytic function $\varphi$ maps the disk onto a convex domain, then $\varphi$ has a continuous extension to the closed unit disk (if $\infty$ is allowed as a value) and assumes no finite value more than once.

\par\smallskip

\subsection{Analytic functions of bounded boundary rotation}

Let us use $V_k$, $k\geq 2$, to denote the set of locally univalent analytic functions  $\varphi$ in the unit disk normalized by the conditions $\varphi(0)=1-\varphi'(0)=0$ that satisfy
\[
\int_{0}^{2\pi}\left|{\rm Re} \left(1+\frac{z\varphi''(z)}{\varphi'(z)}\right)\right|\,d\theta\leq k\pi
\]

for $|z|=|re^{i\theta}|<1$.

\par\smallskip It is known that $\varphi\in V_k$ if and only if

\[
\varphi'(z)={\rm exp}\left\{\frac{1}{\pi}\int_0^{2\pi} \log(1-ze^{-it})^{-1} \, d\mu(t)\right\}\,,
\]
where $\mu(t)$ is a real-valued function  of bounded variation on $[0, 2\pi]$ with
\[
\int_0^{2\pi} d\mu(t)=2\pi \quad\text{and} \quad \int_0^{2\pi} |d\mu(t)|\leq k\pi\,. 
\]

The functions in $V_k$ are precisely those locally univalent analytic functions $\varphi$ normalized as in the class $\mathcal S$ which map $\D$ onto a domain of boundary rotation at most $k\pi$.  We refer the reader to \cite{Pa, Pa2} for the basic properties of the functions in $V_k$. We mention that if $k_1\leq k_2$, then $V_{k_1}\subset V_{k_2}$. We refer the reader to \cite{1}, where analytic functions of bounded boundary rotation are used to generate close-to-convex harmonic mappings.

\par\smallskip
The class $V_2$ consists of functions $\varphi\in\mathcal S$ that are convex. It was proved in \cite{Pa} that if $\varphi\in V_k$ for $2\leq k\leq 4$, then $\varphi$ is univalent.  Pinchuk shows in \cite{Pinchuk} that, indeed, any function in $V_4$ is close-to-convex. 

\par\smallskip
Brannan proves the following (see \cite[Thm. 3.9]{Brannan}). 

\begin{other}\label{thm-Brannan}
Let $\varphi\in V_k$. Assume that $\psi$ is an analytic function in the unit disk with $\psi(0)=\psi'(0)-1=0$ and such that
\begin{equation}\label{eq-ref-Brannan}
\psi'(z)=\varphi'(z)\frac{1+z}{1-z}\,,\quad z\in\D\,.
\end{equation}
Then, $\psi\in V_{k+2}$.
\end{other}

Notice that if $\varphi$ in the previous theorem turns out to be convex (so that $\varphi\in V_2$), the function $\psi$ defined by \eqref{eq-ref-Brannan} belongs to $V_4$. 

\par\smallskip
By adapting the proof of \cite[Thm. 3.9]{Brannan}, we can prove the following more general version of Theorem~\ref{thm-Brannan}.

\begin{cor}\label{cor-Brannan}
Let $\varphi\in V_k$. Assume that $\psi$ is an analytic function in the unit disk with $\psi(0)=\psi'(0)-1=0$ and such that, for some $|\lambda|=1$ and some positive integer $N$, 
\begin{equation*}
\psi'(z)=\varphi'(z)\frac{1-\lambda z^N}{1+\lambda z^N}\,,\quad z\in\D\,.
\end{equation*}
 Then, $\psi\in V_{k+2N}$.
\end{cor}

\begin{pf}
A straightforward calculation gives
\[
1+\frac{z\psi''(z)}{\psi'(z)}=1+\frac{z\varphi''(z)}{\varphi'(z)}-\frac{2\lambda Nz^N}{1-\lambda^2z^{2N}}\,.
\]
Notice that for $z=re^{i\theta}$ and $0\leq r <1$ we have
\begin{align*}
\int_0^{2\pi} & \left|{\rm Re} \left(\frac{2\lambda z^N}{1-\lambda^2z^{2N}}\right)\right|\,d\theta = \frac{1}{2}   \int_0^{2\pi}\left|{\rm Re} \left(\frac{1+\lambda z^N}{1-\lambda z^{N}}-\frac{1-\lambda z^N}{1+\lambda z^{N}} \right)\right|\, d\theta\\
&\leq
 \frac 12 \int_0^{2\pi} {\rm Re}\left(\frac{1+\lambda z^N}{1-\lambda z^{N}}\right)\, d\theta + \frac 12 \int_0^{2\pi} {\rm Re}\left(\frac{1-\lambda z^N}{1+\lambda z^{N}}\right)\, d\theta=2\pi N\,.
\end{align*}
Therefore, we obtain 
\[
\int_{0}^{2\pi}\left|{\rm Re} \left(1+\frac{z\psi''(z)}{\psi'(z)}\right)\right|\,d\theta\leq (k+2N)\pi.
\]
\end{pf}

\par\smallskip

To finish this section, we mention the following theorem that will be fundamental for establishing our main result \cite[Thm. 2.1]{Noonan}.

\begin{other}\label{thm-continuity}
Let $\varphi\in V_k$. Then there exist finite sets 
\[
E_1=\{e^{i\theta_1}, \ldots, e^{i\theta_n} \}
\]
and $E_2\subset E_1$ of points on $|z|=1$, where $n\leq [1+k/2]$, with the following properties: 

\begin{itemize}
\item[(i)] $\varphi(z)$ has continuous (finite-valued) extension to $\overline\D$ except at points of $E_2$; for each $z_j\in E_2$, $|\varphi(z)|\to\infty$ uniformly as $z\to z_j$ in $|z|\leq 1$.

\item[(ii)] $\varphi(z)$ is absolutely continuous on any closed subset of $|z|=1$ disjoint from $E_1$. On any such closed subset, 
\[
\frac{d\varphi(e^{i\theta})}{d\theta}=ie^{i\theta} \varphi'(e^{i\theta})\quad \text{a.e.}\,,
\]
where $d\varphi/d\theta$ is the derivative with respect to values on $|z|=1$ and $\varphi'(e^{i\theta})=\lim_{r\to1^-} \varphi '(re^{i\theta})$.
\end{itemize}
\end{other}

\section{Main result}

The methods used in \cite{HS2}, \cite{B},  and \cite{FHM} to determine whether the vertical shears of a given convex analytic function $\varphi$ normalized by $\varphi(0)=1-\varphi'(0)=0$ are always-convex are based on verifying that the condition in Theorem~\ref{thm-convexcharact} is satisfied. The use of Theorem~\ref{prop-rotations} has allowed us to show that the vertical shears of rotation $\mathcal H_{-1}$, as defined in \eqref{eq-rotation}, of the function $\mathcal H$ given in \eqref{eq-H} are always-convex too. In order to determine \emph{all} convex functions $\varphi$ normalized as above such that the vertical shears of $\varphi$ are always-convex, we need to apply different techniques. In this section, we introduce the relevant previous results in the literature and prove the main theorem in this paper. 
\par\medskip

\subsection{Some useful results on convex harmonic mappings}
Different properties of orientation-preserving harmonic mappings whose dilatation $\omega$ is a rotation of the unit disk, that is, with $\omega(z)=\lambda\, z$ for some $|\lambda|=1$ or, more generally, a finite Blaschke product $B$ of degree $N$ of the form
\begin{equation}\label{eq-Blaschke}
B(z)=e^{i\gamma}\prod_{k=1}^N \frac{a_k-z}{1-\overline{a_k}z}\,,
\end{equation}
where $\gamma\in\R$, and $a_1, a_2, \ldots, a_N\in\D$, have been studied by various authors 
(see, for instance, \cite{HS3, Laugesen, BLW} and the references therein). 

\par\smallskip
It is known that if $f=h+\overline g \in S_H^0$ has constant boundary values, then  $f$ maps the unit disk onto a polygon and its dilatation is a finite Blaschke product. This was proved by Sheil-Small in \cite{S-S}. Hengartner and Schober established the converse. More precisely, they proved the following ``surprising result'' as qualified by themselves in \cite[Thm.~ 3.3]{HS3}.

The \emph{radial  boundary values} $\widehat f$ of a given function $f\in S_H^0$ are defined for $e^{i\theta}\in\partial\D$ by 
\begin{equation}\label{eq-radiallim}
\widehat f(e^{i\theta})=\lim_{r\to 1^-} f(re^{i\theta})\,,
\end{equation}
while the unrestricted limit at such a point $e^{i\theta}$, denoted again by $f$, is equal to
\[
 f(e^{i\theta})=\lim_{z\to e^{i\theta},\ z\in\D} f(z)\,. 
\]

\begin{other}\label{thm-HS}
Let $f$ be a univalent harmonic mapping in the unit disk with dilatation a finite Blaschke product $B$ defined in \eqref{eq-Blaschke}. Assume that $f$ maps $\D$ into a bounded convex domain $D$ and has radial boundary values $\widehat f \in\partial D$ almost everywhere on $\partial \D$. Then, there exists a set $E\subset \partial \D$ with at most $2+N$ points such that the unrestricted limit $f(e^{it})$ exists and is constant on each component of $\partial \D\setminus E$. 
\end{other}

Notice that the requirement that the radial boundary values $\widehat f\in\partial D$ given in \eqref{eq-radiallim} is satisfied under the assumption that $f$ maps the unit disk \emph{onto} a convex domain.
\par\smallskip

The following theorem states some of the results from  \cite[Thm. 2.4]{AM-S}. 

\begin{other}\label{thm-A-MS}
  Let $f=h+\overline g$ be a univalent, harmonic, orientation-preserving mapping from the unit disk $\D$ onto an unbounded convex domain $D$ which is neither a strip nor a half-plane. Then 
\begin{itemize}
    \item[(a)] $f\in h^1$. That is, 
    \[
    \int_{0}^{2\pi} |f(re^{it})|\, dt
    \]
    remains bounded as $r\to 1^-$;
     \item[(b)]There is only one point $e^{i\lambda}$ that corresponds to $\infty$;

     \item[(c)]
     \[
     f(z)=\int_0^{2\pi} P(z,t) \widehat {f}(e^{it})\, dt + A P(z,\lambda)
     \]
     for some constant $A\in\C$. Here, $\widehat f(e^{it})$ denotes the radial limit defined in \eqref{eq-radiallim} and $P(z,\lambda)$ is the Poisson kernel defined by 
     \[
     P(z,\lambda)=\frac{1}{2\pi}{\rm Re\,}\left(\frac{e^{i\lambda}+z}{e^{i\lambda}-z}\right)\,, \quad z \in \D\,, \quad \lambda \in \R \,;
     \]
    \item[(d)] There is a countable set $E\subset \partial\D\setminus\{e^{i\lambda}\}$ such that
    \begin{itemize}
        \item[(i)] The unrestricted limit $\lim f(z)$ exists as $z\to e^{i\theta}$, $z\in\D$, and is continuous for all points $e^{i\theta}\in\partial \D\setminus[E\cup\{e^{i\lambda}\}]$,
        \item[(ii)] The limits 
        \[ \lim_{t\to\theta^-} \widehat f(e^{it})\quad \text{and}\quad \lim_{t\to\theta^+} \widehat f(e^{it})\]
        exist and are different for $e^{i\theta}\in E$.
        \item[(iii)] and the cluster set of $f$ at $e^{i\theta}\in E$ is the line segment joining 
        \[ \lim_{t\to\theta^-} \widehat f(e^{it})\quad \text{to}\quad \lim_{t\to\theta^+} \widehat f(e^{it})\,.
        \]

    \end{itemize}
\end{itemize}
  
\end{other}

It is worth pointing out that the unrestricted limit in item (d) (i) in Theorem~\ref{thm-A-MS} is not necessarily constant at points $e^{i\theta}\in\partial \D\setminus[E\cup\{e^{i\lambda}\}]$. Moreover, in principle, the set $E$ could be dense on the unit circle.  Laugesen \cite{Laugesen} provides further insights into this question in cases when $\partial\D\setminus[E\cup\{e^{i\lambda}\}]$ contains an arc, as demonstrated by the following theorem. 

\begin{other}\label{thm-Laugesen}
Let $I$ be an arc on $|z|=1$. Let $f$ be an orientation-preserving univalent harmonic map in the unit disk such that  $f(e^{it})$ is absolutely continuous on $I$. Then, either $f(e^{it})$ is constant on $I$ or $f$  cannot have dilatation of absolute value $1$ almost everywhere on $I$ while mapping $I$ onto a (strictly) convex arc.
\end{other}
\par\medskip

\subsection{The main theorem} \label{section-main}

We are now ready to state and prove our main result, which is the following theorem. 

\begin{thm}\label{thm-main}
Let $\varphi$ be a convex analytic function in the unit disk with $\varphi(0)=0$ and $\varphi'(0)=1$. Assume that every harmonic mapping  $f=h+\overline g \in S_H^0$ with 
\[
\left\{\begin{array}{cc}
h+g=\varphi\\
g^\prime/h^\prime= \omega
\end{array}\right.\,,
\]
where $\omega$ is a Schwarz function, is convex. Then $f$ is a rotation of a function $F=H+\overline G$ satisfying that $H+G$ belongs to the family $\mathcal A$ defined by \eqref{eq-classA}.
\end{thm}

\begin{pf}
We divide our proof into four different cases. 
\par\smallskip

\textbf{Case 1}. Suppose that for some Schwarz function $\omega$, the harmonic mapping $f$ maps the unit disk $\D$ onto a half-plane $\mathbb H$. Applying a suitable rotation yields a function $F\in S_H^0$, so we may assume that $\mathbb H$ has the form
\[
\mathbb H=\{z\in\C\colon {\rm Re\,} z >-a\}
\]
for some positive real number $a$.
\par\smallskip
It was proved in \cite{Dorff} that the condition $f\in S_H^0$ implies that $a=1/2$ and $H+G$ is equal to the function $\mathcal H$ defined by \eqref{eq-H}. Hence, the result follows in this case.  

\par\smallskip
\textbf{Case 2}. Assume that the conditions in the previous case do not hold but suppose that for some Schwarz function $\omega$, the harmonic mapping $f$ given by the solution to the system in the statement of the theorem with dilatation $\omega$ maps $\D$ onto a strip $\mathbb S$. Applying a suitable rotation yields a function $F\in S_H^0$, so we may assume that $\mathbb S$ is a vertical strip of the form 
\[
\mathbb S=\{z\in\C\colon a<{\rm Re\,} z <b\}\,, 
\]
where $a$ and $b$ are real numbers with $a<0<b$ and $|a|<b$. Then, as is shown in \cite{GN}, we can write, for some $A>0$ and some $\gamma\in[0, \pi/2)$,
\par\smallskip

\[
\mathbb S=\left\{ z\in\C\colon A\left(\gamma -\frac{\pi}{2}\right)<{\rm Re\,} z <A\left(\gamma +\frac{\pi}{2}\right)\right\}\,.
\]

\par\smallskip
It is a consequence of the results in  \cite{GN} that, within these terms,  
\[
H(z)+G(z)=-iA\log\frac{1+i e^{i\gamma}z}{1-i e^{-i\gamma}z} =-i A \log\frac{1-\overline\lambda z}{1-\lambda z} \,,
\]
where $\lambda=ie^{-i\gamma}\not\in\{-1, 1\}$ and $z\in \D$.

\par\smallskip
The condition $F\in S_H^0$ and the fact that the dilatation of $F$ is a Schwarz function imply $(H+G)'(0)=1$, which gives $A=1/(2{\rm Im\,} \lambda)$. That is, $H+G$ is
equal to one of the functions $\mathcal L_\lambda$ defined by \eqref{eq-llambda}, a case also considered in \cite{Dorff}. Hence, the result holds true in this case as well. 

\par\smallskip
\textbf{Case 3}. Let us assume now that the requirements in the previous cases are not satisfied for any vertical shear of the analytic function $\varphi$. 

Recall that if $\varphi$ maps the unit disk onto a half-plane or onto a strip, then, according to our exposition in Section~\ref{sec-convexanalytic}, $\varphi\in\mathcal A$ and the analysis is complete. Hence, if we assume that $\varphi$ does not map the unit disk onto a half-plane or a strip, only two possibilities may occur:

- If $\varphi$  is unbounded, then there is only one point $e^{i\lambda}$, provided by Theorem~\ref{thm-A-MS}, such that $\varphi(e^{i\lambda})=\infty$. Since $\varphi$ is convex, it follows that $\varphi(e^{it})$, and hence ${\rm Re\,}\varphi(e^{it})$, is continuous at any $t\neq \lambda$. 

- If $\varphi$ is bounded, then ${\rm Re\,}\varphi(e^{it})$ is continuous for all $t\in[0, 2\pi]$.

\par\smallskip
The case we consider at this point is the following: Suppose that the function $f=h+\overline {g}\in\mathcal S_H^0$ obtained from the equation $h+g=\varphi$ and with dilatation $\omega=z$ is bounded. Notice that this assumption does not necessarily imply that $\varphi$ maps the unit disk onto a bounded convex domain.

\par\smallskip
Since $f$ must be convex, we can apply Theorem~\ref{thm-HS} to fix the set $E$ with at most $3$ points such that the unrestricted limit $f(e^{it})$ exists and is constant on each component of $\partial\mathbb{D}\setminus E$.

\par\smallskip
\textbf{Case 3.1.} Suppose that $E$ is empty. In this case, $f(e^{it})$ is constant, say $x_1+iy_1$, on $\partial \mathbb{D}$.

\par\smallskip
Since $f(e^{it})$ exists and is constant for all $|z|=1,$ ${\rm Re\,} f(e^{it})$ exists and is also constant on $\partial \mathbb{D}.$ This implies that ${\rm Re\,} \varphi(e^{it})$ exists and is constant for all $z\in\partial\mathbb{D}$. Therefore, the boundary values of $\varphi$ are contained in a straight vertical line. Since $\varphi(\mathbb{D})$ is convex, we necessarily have that the boundary values \emph{are} the whole vertical line. That is,  $\varphi(\mathbb{D})$ is a half-plane. According to our discussion in Section~\ref{sec-convexanalytic}, $\varphi\in\mathcal A$.

\par\smallskip
\textbf{Case 3.2.} Let $E=\{e^{i\theta_1}\}$ for some $0\leq \theta_1 < 2\pi$. In this case, $f(e^{it})$ is equal to a constant, say $x_1+iy_1$ (which is not necessarily equal to the constant in the previous case), on $[0, 2\pi)\setminus \{\theta_1\}$. As before, this implies that ${\rm Re\,}\varphi(e^{it})=x_1$ for all such $t$. The continuity of the boundary values of convex analytic functions again gives that $\varphi$ maps the unit disk onto a half-plane, so that $\varphi\in\mathcal A$.

\par\smallskip
\textbf{Case 3.3.} We now jointly consider the cases when either $E=\{e^{i\theta_1}, e^{i\theta_2}\}$ for some $0\leq \theta_1< \theta_2 < 2\pi$ or $E=\{e^{i\theta_1}, e^{i\theta_2}, e^{i\theta_3}\}$ for some $0\leq \theta_1< \theta_2 < \theta_3 < 2\pi$. The argument is similar in both cases and the conclusion is obtained as a consequence of the fact that  ${\rm Re\,}\varphi(e^{it})$ turns out to be constant at all points on the unit circle except perhaps one, so that the result follows as before. In fact, the case when $E$ consists of $2$ points can be reduced to the one with $3$ points by setting $e^{i\theta_2}=e^{i\theta_3}$.

Let us then suppose that $f(e^{it})$ is equal to a constant, say $x_1+iy_1$ on $I_1=(\theta_1, \theta_2)$, $x_2+iy_2$ on $I_2=(\theta_2, \theta_3)$, and $x_3+iy_3$  on $(0,2\pi)\setminus[\overline{I_1}\cup \overline{I_2}]$.

This gives that the boundary values ${\rm Re\,} \varphi(e^{it})={\rm Re\,} f(e^{it})$ are given by
\[
{\rm Re\,} \varphi(e^{it})=\begin{cases}
    x_1\,,\quad t\in I_1\,\\
    x_2\,,\quad t\in I_2\,\\
    x_3\,,\quad t\in I_3\,.
    \end{cases}
\]

Using again that ${\rm Re\,}\varphi(e^{it})$ is continuous at all points on the unit circle except perhaps for one gives $x_1=x_2=x_3$, which allows us to conclude, as announced, that  $\varphi\in\mathcal A$.

Before continuing, it is convenient to notice that the choice of $\omega(z)=z$ in Case 3 could have been changed by any other dilatation of the form $\omega(z)=\lambda z^N$, where $|\lambda|=1$ and $N$ is a positive integer.

\par\smallskip
\textbf{Case 4.} The cases that remain to be analyzed are precisely those when all the vertical shears of $\varphi$ with dilatation $\omega(z)=\lambda z,$ where $|\lambda|=1,$ are unbounded and none of them map the unit disk onto a half-plane or a strip.  
\par\smallskip
The key tools used to resolve this case are Theorems ~\ref{thm-Brannan}, \ref{thm-continuity}, and \ref{thm-Laugesen}.
\par\smallskip
More concretely, we are to prove that the vertical shear of $\varphi$ with dilatation $\omega(z)=-z$ satisfies the following condition: There exists a finite set
\[
E=\{e^{i\theta_1}, e^{i\theta_2}, \ldots, e^{i\theta_n}\}
\]
with at most five points such that $f(e^{it})$ is absolutely continuous on any closed subset of $|z|=1$ disjoint from $E$.
\par\smallskip
This will imply, by Theorem~\ref{thm-Laugesen}, that $f(e^{it})$ must be constant in any closed subset of $|z|=1$ disjoint from $E$. With this information at hand, it suffices to argue as in Case 3 to get the desired conclusion. 
\par\smallskip
In order to prove the existence of the set $E$ as above, we make use of Theorem~\ref{thm-continuity} applied to the function $h+g=\varphi\in V_2$. This produces a set $S_1=\{e^{i\theta_1}, e^{i\theta_2}\}$ with at most two points such that $\varphi(e^{it})$ is absolutely continuous on any closed subset of $|z|=1$ disjoint from $S_1$.
\par\smallskip
Since 
\[
h'(z)-g'(z)=\varphi'(z)\frac{1+z}{1-z}\,, z\in\mathbb D\,,
\]
we have, by Theorem~\ref{thm-Brannan}, that $h-g\in V_4$. Therefore, there is a set 
$S_2=\{e^{i\theta_3}, e^{i\theta_4}, e^{i\theta_5}\}$, with at most $3$ points,
such that $(h-g)(e^{it})$ is absolutely continuous on
any closed subset of $|z| = 1$ disjoint from $S_2$. Let $E=S_1\cup S_2$, which has, at most, five points.
\par\smallskip
Since the conjugate of an absolutely continuous function and the sum of absolutely continuous functions are again absolutely continuous, we have that 
\[
h=\frac{\varphi+h-g}{2}\quad \text{and} \quad \overline g=\frac{\overline\varphi-\overline{(h-g)}}{2}
\]
satisfy that $h(e^{it})$ and $\overline{g(e^{it})}$, and consequently $f(e^{it})$, are absolutely continuous on any closed subset of $|z| = 1$ disjoint from $E$. 
\par\smallskip
This proves our claim and ends the proof of the theorem. 
\end{pf}

\par\smallskip
We conclude with the following remark. Notice that the proof of Theorem~\ref{thm-main} can be adapted to those cases when $\varphi$ is convex and the vertical shear of $\varphi$ with dilatation $\omega(z)=\lambda z^N$, for some positive integer $N$ and some $|\lambda|=1$ is convex. Hence, our initial goal of characterizing the shears of $\varphi$ in the $\theta$ direction that remain convex for every Schwarz function $\omega$ leads to an important consequence: the characterization of the shears of a convex function $\varphi$ with dilatation $\omega$ of the form mentioned, which turned out to be convex. That is, we have the following corollary.

\begin{cor}
Let $\varphi$ be a convex analytic function in the unit disk with $\varphi(0)=0$ and $\varphi'(0)=1$. Let $N$ be a positive integer and let $|\lambda|=|\eta|=1$. Assume that the harmonic mapping  $f=h+\overline g \in S_H^0$, defined in the unit disk by the system of equations  
\[
\left\{\begin{array}{cc}
h-\eta g=\varphi\\
g^\prime(z)/h^\prime(z)= \omega(z)=\lambda z^N
\end{array}\right.\,,
\]
is convex.  Then, $f$ is a rotation of a function $F=H+\overline G$ satisfying that $H+G$ belongs to the family $\mathcal A$ defined in \eqref{eq-classA}.
\end{cor}

\begin{pf}
By applying a rotation, we can assume, without loss of generality, that $\eta=-1$. In order to prove the corollary, it suffices to check that the set $E$ in Case $4$ in the proof of Theorem~\ref{thm-main} is finite. This condition will be derived from Theorem~\ref{thm-continuity} if it is proved that $h-g \in V_k$ for some $k\geq 2$.
Since
\[
h'(z)-g'(z)=\varphi'(z) \left(\frac{1-\lambda z^N}{1+\lambda z^N}\right)
\]
and $\varphi\in V_2$, the result follows from Corollary~\ref{cor-Brannan}.
\end{pf}

\section*{Acknowledgements}

M. Wo{\l}oszkiewicz-Cyll would like to thank the Department of Mathematical Analysis at the University of La Laguna for their hospitality during her research visit. 

The authors would like to thank the referees for their careful reading of the manuscript, for their useful suggestions to improve the exposition, and for providing fundamental references for the bibliography. They are also grateful to Professors M. D. Contreras and L. Rodr\'iguez-Piazza for bringing to their attention an inaccuracy in the first paragraph of an earlier version of the manuscript, which also appears in the published version but does not affect the results. This has been corrected in the current arXiv submission.

\end{document}